# A new derivation of Hermite's integral for the Hurwitz zeta function

Donal F. Connon

dconnon@btopenworld.com

11 August 2009

**Abstract**

We obtain another proof of Hermite's integral for the Hurwitz zeta function $\varsigma(s,u)$.

**Proof**

Chen [5] has recently reported that for $u > 0$ and $s > 0$

(1) $$\frac{2}{\Gamma(s)}\int_0^\infty e^{-uy^2} y^{2s-1} \sin(xy^2)\,dy = \frac{\sin[s\tan^{-1}(x/u)]}{(u^2+x^2)^{s/2}}$$

and multiplying this across by $\dfrac{1}{e^{2\pi x}-1}$ and integrating with respect to $x$ results in

(2) $$\frac{2}{\Gamma(s)}\int_0^\infty \frac{1}{e^{2\pi x}-1}dx \int_0^\infty e^{-uy^2} y^{2s-1} \sin(xy^2)\,dy = \int_0^\infty \frac{\sin[s\tan^{-1}(x/u)]}{(u^2+x^2)^{s/2}(e^{2\pi x}-1)}\,dx$$

Reversing the order of integration we note that

$$\int_0^\infty \frac{1}{e^{2\pi x}-1}dx \int_0^\infty e^{-uy^2} y^{2s-1} \sin(xy^2)\,dy = \int_0^\infty e^{-uy^2} y^{2s-1}\,dy \int_0^\infty \frac{\sin(xy^2)}{e^{2\pi x}-1}\,dx$$

and using Legendre's relation [7, p.122]

$$2\int_0^\infty \frac{\sin(xt)}{e^{2\pi x}-1}\,dx = \frac{1}{e^t-1} - \frac{1}{t} + \frac{1}{2} = \frac{1}{2}\coth\frac{t}{2} - \frac{1}{t}$$

(a rigorous derivation of this result is shown in Bromwich's book [4, p.501]), we obtain

$$\int_0^\infty \frac{1}{e^{2\pi x}-1}dx \int_0^\infty e^{-uy^2} y^{2s-1} \sin(xy^2)\,dy = \frac{1}{2}\int_0^\infty e^{-uy^2} y^{2s-1}\left[\frac{1}{e^{y^2}-1} - \frac{1}{y^2} + \frac{1}{2}\right]dy$$

With the substitution $v = y^2$ this becomes

$$= \frac{1}{4}\int_0^\infty e^{-uv} v^{s-1}\left[\frac{1}{e^v-1}-\frac{1}{v}+\frac{1}{2}\right]dv$$

We now recall the well-known formula for the Hurwitz zeta function which is reported in [6, p.92] as being valid for $\text{Re}(s) > -1$

(3) $$\varsigma(s,u) = \frac{u^{-s}}{2} + \frac{u^{1-s}}{s-1} + \frac{1}{\Gamma(s)}\int_0^\infty e^{-uv} v^{s-1}\left[\frac{1}{e^v-1}-\frac{1}{v}+\frac{1}{2}\right]dv$$

and we thereby obtain Hermite's integral [1, p.55]

(4) $$\varsigma(s,u) = \frac{u^{-s}}{2} + \frac{u^{1-s}}{s-1} + 2\int_0^\infty \frac{\sin[s\tan^{-1}(x/u)]}{(u^2+x^2)^{s/2}(e^{2\pi x}-1)}dx$$

□

Chen [5] has stated that (1) is valid for $s > 0$ and we now consider the limit as $s \to 0$. We easily see that

$$\Gamma(s)\sin[s\tan^{-1}(x/u)] = \frac{s\tan^{-1}(x/u)\Gamma(s)\sin[s\tan^{-1}(x/u)]}{s\tan^{-1}(x/u)}$$

$$= \tan^{-1}(x/u)\Gamma(s+1)\frac{\sin[s\tan^{-1}(x/u)]}{s\tan^{-1}(x/u)}$$

and hence we have

$$\lim_{s \to 0}\Gamma(s)\sin[s\tan^{-1}(x/u)] = \tan^{-1}(x/u)$$

Therefore we deduce that

$$2\int_0^\infty \frac{e^{-uy^2}\sin(xy^2)}{y}dy = \tan^{-1}(x/u)$$

or equivalently we obtain the well-known integral

(5) $$\tan^{-1}(x/u) = \int_0^\infty \frac{e^{-uy}\sin(xy)}{y}dy$$

Rigorous derivations of (5) are contained in [2, p.285] and [3, p.272].

**REFERENCES**

[1] G.E. Andrews, R. Askey and R. Roy, Special Functions.
    Cambridge University Press, Cambridge, 1999.

[2] T.M. Apostol, Mathematical Analysis, Second Ed., Addison-Wesley Publishing Company, Menlo Park (California), London and Don Mills (Ontario), 1974.

[3] R.G. Bartle, The Elements of Real Analysis. 2$^{nd}$ Ed. John Wiley & Sons Inc, New York, 1976.

[4] T.J.I'a Bromwich, Introduction to the theory of infinite series. 2$^{nd}$ edition Macmillan & Co Ltd, 1965.

[5] H. Chen, The Fresnel integrals revisited.
    Coll. Math. Journal, 40, 259-262, 2009.

[6] H.M. Srivastava and J. Choi, Series Associated with the Zeta and Related Functions. Kluwer Academic Publishers, Dordrecht, the Netherlands, 2001.

[7] E.T. Whittaker and G.N. Watson, A Course of Modern Analysis: An Introduction to the General Theory of Infinite Processes and of Analytic Functions; With an Account of the Principal Transcendental Functions. Fourth Ed., Cambridge University Press, Cambridge, London and New York, 1963.

Donal F. Connon
Elmhurst
Dundle Road
Matfield
Kent TN12 7HD
dconnon@btopenworld.com